# Predictable trajectories of the reduced Collatz iteration and a possible pathway to the proof of the Collatz conjecture (Version 2).

Leonel Sternberg

## Abstract

I show here that there are three different kinds of iterations for the reduced Collatz algorithm; depending on whether the root of the number is odd or even. There is only one kind of iteration if the root is odd and two kinds if the root is even. I also show that iterations on numbers with odd roots will cause an increase in value and eventually lead to an even rooted number. The iterations on even rooted numbers will subsequently cause a decrease in values. Because increase in values during the odd root iterations are bounded, I conclude that the Collatz iteration cannot veer to infinity. Since the sequence generated by the Collatz iteration is infinite and the values of the numbers do not veer to infinity it must either cycle and/or converge. I postulate that any cycling must occur with alternating types of iterations: e.g. odd rooted iterations which cause the values of the numbers to increase followed by even rooted iterations which causes the values to decrease. I show here that for simpler types of cycles, valid values of odd rooted or even rooted numbers are only found in a narrow gap which closes as the number of iterations increase. I further generalize to all types of odd-even and even-odd iterations. Given that previous work has shown that only very large non-trivial cycles are feasible during the Collatz iteration and this study shows the low probability of large simple cycles, leads us to the conclusion most likely cycles other than the trivial cycle are not possible during the Collatz iteration.

## Introduction

Hailstone sequences are generated by the Collatz iteration which specifies the following algorithm

$$C(x_{i+1}) \quad \begin{cases} = \dfrac{x_i}{2} \, for \; x \equiv 0 \; (mod \; 2) \\ = 3x_i + 1 \; for \; x \; \equiv 1 \; (mod \; 2) \end{cases}$$

The sequence of numbers generated by this iteration often seems chaotic with wild swings and unpredictable number of iterations to reach the value of 1. But the Collatz conjecture states that eventually regardless of the initial number the sequence will decrease to 1 and cycle between 1 and 4. This conjecture has, as yet, not been proven. Computer calculations, however, have not been able to find any exceptions to this conjecture (Lagarias 2010). There are other variations of the Collatz iteration. Here, I use the reduced Collatz iteration (Colussi 2011):

$$R(x_{i+1}) = \frac{3x_i + 1}{2^a},$$

by which $a$ is the highest power of 2 that will divide $3x + 1$. The sequence generated by the reduced Collatz iteration is composed only of odd numbers, but it shows the same chaotic fluctuation as the Collatz iterated sequence. Note that with the reduced Collatz iteration there is no trivial cycle between 1 and 4, as even numbers are not included in the iteration output. Rather, the trivial cycle as expressed by the reduced Collatz iteration is an identity function or a convergence to 1:

$$R(1) = \frac{3 \cdot 1 + 1}{2^2} = 1$$

There are several questions regarding the above sequence, which are peripheral to an actual proof of the conjecture, that need answers. For example, are there cycles other than the trivial (1-4) cycle? Why does the sequence fluctuate wildly? Starting with the value of 27 the sequence can be as high as 3077 before converging to 1 after 41 iterations (Table 1). Why does the sequence eventually converge to 1? Here, I show that there are three major predictable pathways of the hailstone sequence. In some case, one can predict the value of the number after several iterations without actually going through the above algorithm calculation; i.e. multiplying the number by 3 adding one and searching for the highest power of 2 that will divide the number. Further, I show that the probability of large cycles other than the trivial cycle is miniscule during the Collatz iteration. The work presented here suggests a possible pathway to prove the Collatz conjecture.

**The trajectory of the hailstone sequence depends on the root of the number.**

*Definition:* Given any odd number (*m*) we define the root of this number as *a* where

$$m = 2a + 1.$$

Now *a*, the root, can be either even or odd and this will affect the type of reduced Collatz iteration carried out on the number. Numbers with odd roots will cause the hailstone sequence to increase after each iteration, while numbers with even roots will cause the hailstone sequence to decrease after each iteration. Further these iterations are predictable.

**Theorem 1.** Given an odd number *m*, if its root *a* is odd, then the next value of the sequence will be greater than that *m* after one reduced Collatz iteration; i.e.

$$\frac{3 \cdot m + 1}{2^x} > m\,,$$

$$if\ m\ has\ an\ odd\ root\ and\ x\ is\ the\ largest\ exponent\ of\ 2\ that\ divides\ the\ numerator$$

Proof:

Let $m = 2a + 1$ where $a = 2b + 1$, i.e. the root is odd

Therefore $m = 2^2b + 3$

Carrying out the Collatz iteration:

$$\frac{3 \cdot (2^2 b + 3) + 1}{2^x}$$

Where *x* is the maximum exponent of 2 dividing the numerator. The above expression is simplified to:

$$\frac{3 \cdot 2^2 b + 2 \cdot 5}{2^x} = 3 \cdot 2 \cdot b + 5 > 2^2 b + 3 = m$$

Q.E.D.

On the other hand, if the root of the odd number is even, the Collatz iteration will produce a lower value compared to the original number.

**Theorem 2.** Given an odd number *m*, if it's root is even, then the next value of the sequence will be less than *m* after one reduced Collatz iteration; i.e.

$$\frac{3 \cdot m + 1}{2^x} < m$$

Let $m = 2a + 1$ where $a = 2b$ for some odd number b

Therefore $m = 2^2 b + 1$

Proof:

$$\frac{3 \cdot (2^2 b + 1) + 1}{2^x} = \frac{3 \cdot 2^2 \cdot b + 4}{2^x} = \frac{3b + 1}{2^{x-2}} < 2^2 b + 1 = m$$

Q.E.D.

In conclusion: for any odd number (*m*) after the reduced Collatz iteration, the number in the sequence will either be > *m* if the root of *m* is odd and <*m* if the root of *m* is even.

The effect of several iterations on the root is very predictable for all cases with the exception of one particular case. By predictable, I mean that the new value of the root is predicted without actually going through the Collatz iteration algorithm. The next theorems demonstrate how the iterations of odd and some types of even roots are predictable.

*Definition:* Given either an even or odd root ($2^y \cdot 3^u \cdot b$ or $2^y \cdot 3^u \cdot b - 1$, $where\ u \geq 0$ respectively) with *b* odd not a multiple of 3. The factor $2^y \cdot 3^u$ is called here the multiplicative factor of the root.

**Theorem 3.** If *m* is an odd number with an odd root, then the multiplicative factor of the odd root increases by $3 \cdot 2^{-1}$ after each iteration.

Proof:

Let *a* be an odd root $= 2^y \cdot 3^u \cdot b - 1$ with $b$ an odd integer non multiple of 3 and y≥1. Note that the odd number a is expressed somewhat different than the standard *2n+1* because it is easier to do the iterations with this format. Our format is simply *2(n+1)-1*.

$m = 2(2^y \cdot 3^u \cdot b - 1)$+1=$2^{y+1} \cdot 3^u \cdot b$ −1

By iterating:

$$\frac{3(2^{y+1} \cdot 3^u \cdot b - 1) + 1}{2^x}$$

Where *x* is the largest exponent of 2 that will divide the numerator.

The above equation becomes:

$$\frac{2^{y+1} \cdot 3^{u+1} \cdot b - 2}{2^x} = 2^y \cdot 3^{u+1} \cdot b - 1 = 2(3^{u+1} \cdot 2^{y-1} \cdot b - 1) + 1$$

With the new root $(3^{u+1} \cdot 2^{y-1} \cdot b - 1)$ having the multiplicative factor 3/2 greater than the previous root.

Q.E.D.

Note that the root $(3^{u+1} \cdot 2^{y-1} \cdot b - 1)$ is still odd if *y-1≥1*. Further iterations will occur until an even root is derived:

$$3^{u+y} \cdot b - 1$$

In short: a number with an odd root will eventually beget a number with an even root after a finite set of iterations which is the subject of the next theory. Further, each iteration of an odd rooted number will increase the value of the root multiplicative factor by 1.5 until it reaches an even root value.

**Theorem 4.** Given an odd number $m$ with an odd root $2^y \cdot 3^u \cdot b - 1$. After *y* iterations, the derived odd number will have an even root $3^{u+y} \cdot b - 1$.

This follows logically from the comments above.

An odd number having odd root results in an increasing sequence of numbers after each iteration. At some point, it will reach an even root, and subsequently the sequence will start to decrease. Number 31=2*15+1 has a root $2^4 \cdot 1 - 1$. Therefore, there will be four iterations and the result will be an even root $3^4 \cdot 1 - 1 = 80$ for the number 161generated 4 iterations after the number 31.

**Theorem 5.** If *m* is an odd number with an even root $= 2^u \cdot b$. With b an odd number and *u>2*, then the multiplicative factor of the root after one iteration will decrease by $3 \cdot 2^{-2}$, but the root remains even.

Proof: let *a* be an even root $= 2^u b$, with $b$ and odd number and *u>2*. Therefore

$$m = 2(2^u b) + 1 = 2^{u+1} b + 1$$

The iteration of *m* is:

$$\frac{3(2^{u+1}b + 1) + 1}{2^x} = \frac{3 \cdot 2^{u+1} \cdot b + 4}{2^x} = 3 \cdot 2^{u-1} \cdot b + 1 = 2(3 \cdot 2^{u-2} \cdot b) + 1$$

The new root $(3 \cdot 2^{u-2} \cdot b)$ is even as long as *u>2*.

Q.E.D.

The iterations can continue, but the outcome whether the exponent *u* is an even or odd will differ. If *u* is an even number eventually after *u/2* iterations the root will become $3^{u/2}b$ which is an odd root. However, if *u* is an odd number after *(u-1)/2* iterations the exponent of 2 will be 1. Which has a variable outcome with regards to the next root.

**Theorem 6**. If *m* is an odd number with and even root *=2b* where *b* is odd. Then the value of the new root will be:

$$\frac{3b + 1 - 2^{x-2}}{2^{x-1}}, \qquad x \geq 3$$

Proof: *m=2(2b)+1*

Carrying out the iteration

$$\frac{3(2^2 b + 1) + 1}{2^x} = \frac{3 \cdot 2^2 \cdot b + 4}{2^x} = \frac{3b + 1}{2^{x-2}}$$

One can see that $x \geq 3$

The new root is calculated as:

$$\frac{3b + 1}{2^{x-2}} = 2a + 1 \rightarrow a = \frac{3b + 1 - 2^{x-2}}{2^{x-1}}$$

A general scheme for all the possible results of the reduced iterations in the hailstone sequence with different types of roots is shown in Figure 1. Of all the possibilities, there is only one with an unpredictable outcome (i.e. one actually has to go through the reduced Collatz iteration to calculate changes in the root value, labeled as type II even iteration). The unpredictable outcome is given by odd numbers having an even root *=2b* where *b* is an odd number. All other outcomes are predicted by Theorems 4 and 5. Theorem 4 also shows that any increase in values of the numbers in the hailstone sequence has an upper boundary. Collatz iterations on the odd number having the root $2^y \cdot 3^u \cdot b - 1$ will cause the root (and subsequently the number) to increase in value but the root upper boundary is $3^{u+y} \cdot b - 1$. This is important since it leads to the conclusion that the hailstone sequence cannot veer off to infinity, as any increase is bounded.

The hailstone sequence is an infinite sequence and since no value in the sequence is infinite it must converge and or cycle. The observation that odd rooted number will increase after iterations and even rooted numbers will decrease, is also important on determining whether the Hailstone sequence could cycle. Any cycling must have an increase in the value of the number followed by a decrease in the value; i.e. odd rooted iterations followed by even rooted iterations. On the next section we will examine the possibility of cycling using odd rooted iterations followed by even rooted iterations. We do it in this order but we can show the lack of cycling just as well with even rooted iterations followed by odd rooted iterations.

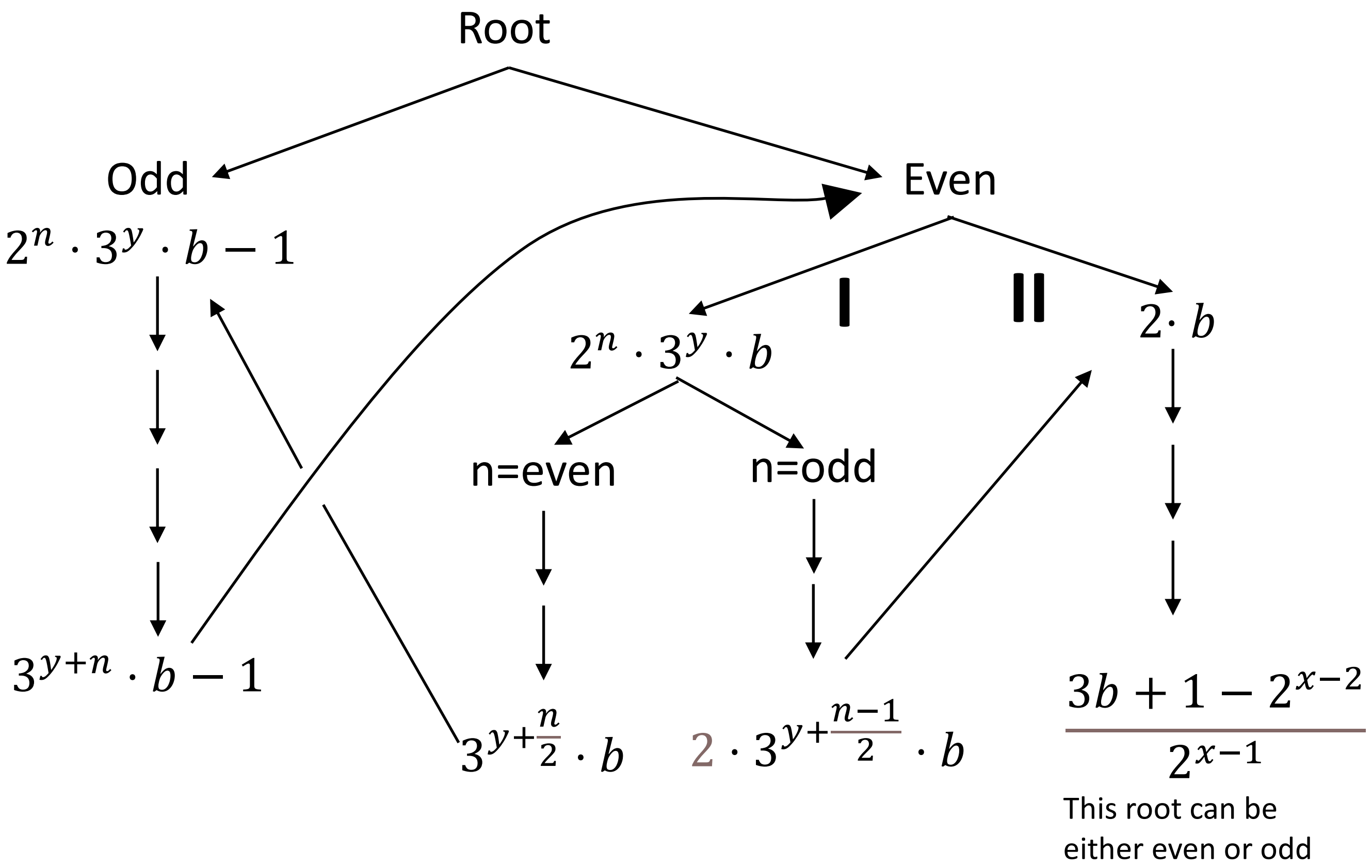

**Figure 1.** Flow chart of possible outcome of the reduced Collatz iterations of odd numbers with either odd or even roots.

**On the possibility of cycling in hailstone sequence**.

Given the flow chart of trajectories (Figure 1), we can conclude that any cycle has to involve iterations off one type of root (for example odd rooted numbers) followed by iterations of the other type of roots (even in the case where one starts with an odd rooted number). No cycle can occur with only odd root iterations as the number increases progressively with each iteration. Likewise, there can be no cycle with only even root iterations as the number will progressively decrease with each even root iteration. Here I will delineate three of possible combinations of odd-even trajectories and show the limits of where potential cycles can occur during the reduced Collatz iteration.

**First case**:

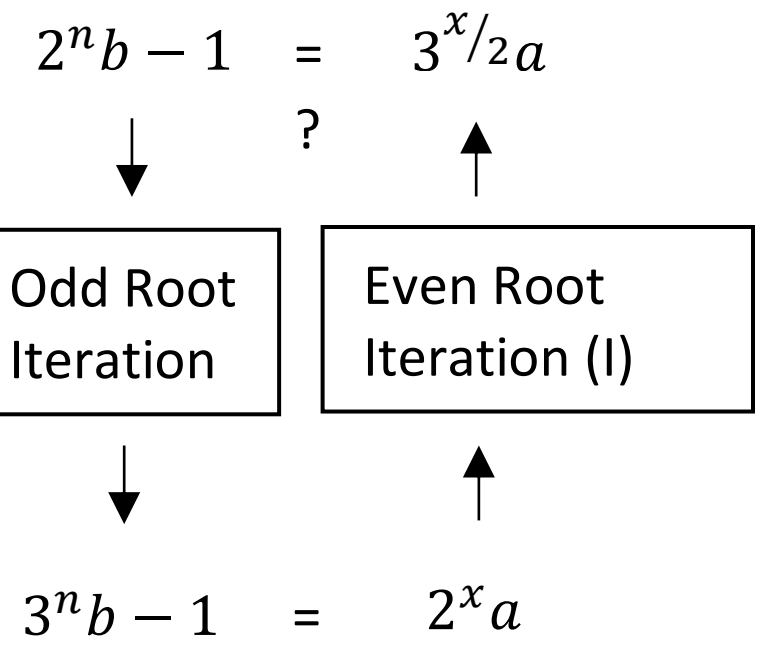

The first case to be considered is to start with an odd rooted number followed by odd root iterations until it's even rooted (with x>2 and even), followed by even rooted iteration type I until its odd rooted again. Will the initial odd root ($2^n b - 1$) ever be equal to the newly generated odd root ($3^{x/2} a$ )? I will show that, if this were the case, there are limits on the values of x, which would give a potential value of b greater than 0 and ≥1. Further, I will show that as the number of odd root iterations a gap of potential solutions decreases. Assume that indeed:

$$2^n b - 1 = 3^{x/2} a$$

From the bottom equality in the diagram above we can say that

$$a = \frac{3^n b - 1}{2^x}$$

Substituting the value of a in the equality above:

$$2^n b - 1 = \frac{3^{\frac{x}{2}+n} b - 3^{x/2}}{2^x}$$

Solving for b:

$$b = \frac{2^x - 3^{x/2}}{2^{x+n} - 3^{\frac{x}{2}+n}}$$

I will examine conditions when b is either <0 , <1, which is contrary to the basic definition of an odd positive root ($2^n b - 1$) with b as a positive integer.

<u>First I calculate the values of x which will lead to b <0.</u>

Define $\varphi = log_2 3$

The numerator of the above equation is positive since $2^{x} > 3^{x/2} = 2^{\varphi x/2}$. The denominator therefore, has to be negative: i.e.

$$2^{x+n} < 3^{\frac{x}{2}+n}$$

$$\Rightarrow 2^{x+n} < 2^{\varphi(\frac{x}{2}+n)}$$

$$\Rightarrow x+n < \varphi(\frac{x}{2}+n)$$

And

$$x < \frac{(\varphi-1)n}{(1-\frac{\varphi}{2})} = \frac{log_2 \frac{3^n}{2^n}}{(1-\frac{\varphi}{2})}$$

The right side of the above inequality is in $\log_2$ form to compare it with the equation to be developed below.

<u>Next I calculate the values of x where b<1.</u>

b<1 when the denominator of the equation defining b is greater than the numerator:

$$2^{x+n} - 3^{\frac{x}{2}+n} > 2^{x} - 3^{x/2}$$

$$\Rightarrow 2^{x+n} - 2^{x} > 2^{\varphi\left(\frac{x}{2}+n\right)} - 2^{\varphi x/2}$$

$$\Rightarrow 2^{x}(2^{n}-1) > 2^{\frac{\varphi x}{2}}(2^{\varphi n}-1)$$

$$\Rightarrow \frac{2^{x}}{2^{\frac{\varphi x}{2}}} = 2^{x-\frac{\varphi x}{2}} > \frac{2^{\varphi n}-1}{2^{n}-1}$$

$$\Rightarrow x - \frac{\varphi}{2}x > log_2(\frac{2^{\varphi n}-1}{2^{n}-1})$$

$$\Rightarrow x > \frac{log_2(\frac{3^{n}-1}{2^{n}-1})}{1-\varphi/2}$$

In Summary b will have values contradictory to the premise that it’s a positive odd integer if

$$x < \frac{log_2 \frac{3^n}{2^n}}{(1-\frac{\varphi}{2})},$$

which would cause b <0 or

$$x > \frac{log_2\left(\frac{3^n-1}{2^n-1}\right)}{(1-\frac{\varphi}{2})},$$

which would cause b <1.

The limits of x leading to illegitimate values of b is visualized by the following number line where I stands for integers:

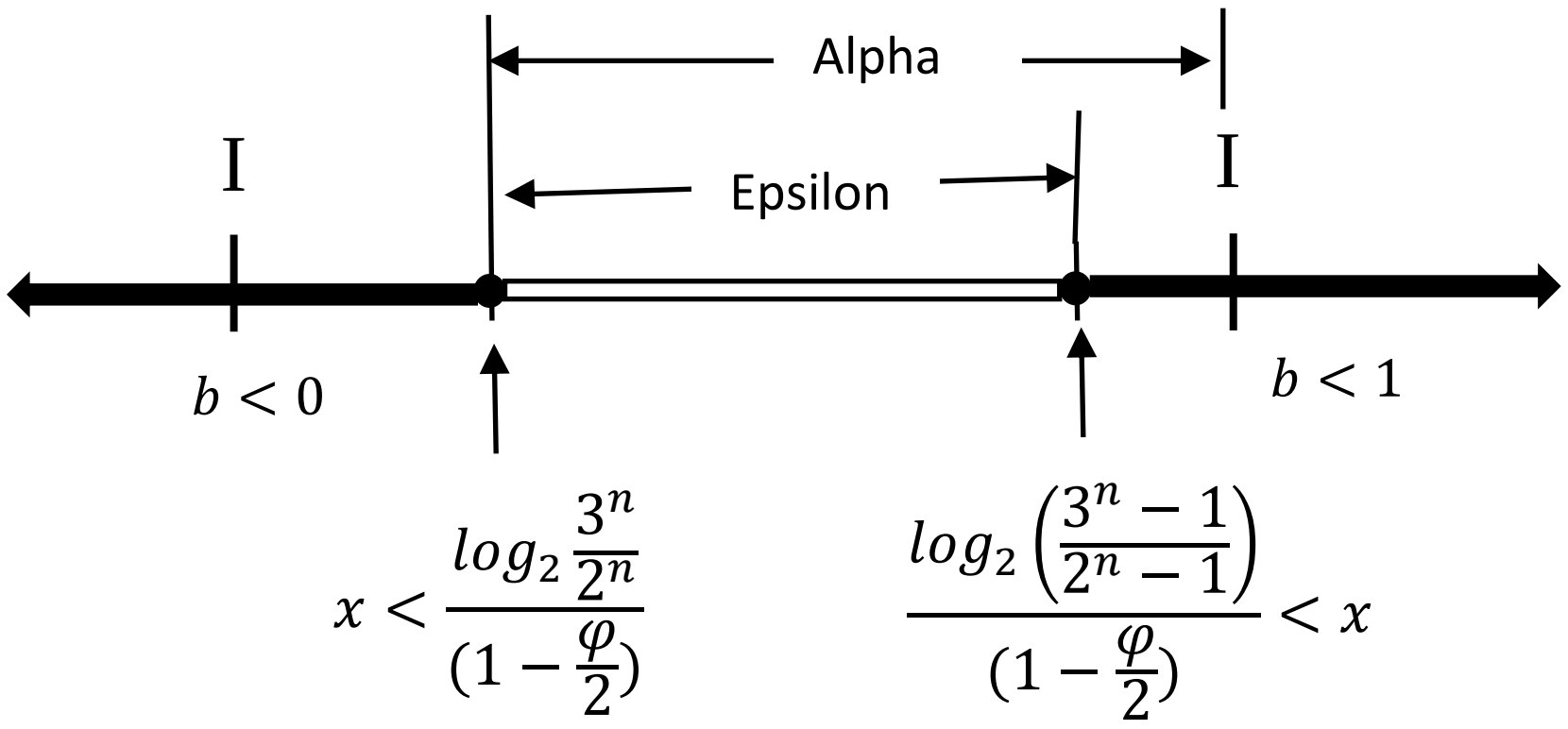

Is easy to see that as n increases the two limit points converge. It must be shown that there are no integers between these limit points i.e. Epsilon<Alpha. Obviously, at n=infinity the gap is completely closed and all values of x's will lead to a contradiction to the basic premise that b is positive integer and no cycles are possible. I calculate Alpha and Epsilon as:

$$Alpha1_{min}(n) = \ min_1^n(1 - mod_1\frac{log_2\frac{3^n}{2^n}}{\left(1-\frac{\varphi}{2}\right)})$$

$$Epsilon1(n) = \frac{log_2\left(\frac{3^n-1}{2^n-1}\right)}{(1-\frac{\varphi}{2})} - \frac{log_2\frac{3^n}{2^n}}{(1-\frac{\varphi}{2})}$$

The Alpha function cycles and so as to better visualize I modify it to $Alpha_{min}$ as shown above. This is a more stringent requirement that $Epsilon(n) < Alpha(n)$. So that satisfying the above equation satisfies $Epsilon(n) < Alpha(n)$. It can be seen that $Epsilon(n) > Alpha_{min}(n)$ for $n \leq 6$, but $Epsilon(n) < Alpha_{min}(n)$ for $n > 6$ (Figure 2). We can manually check for the values of $n \leq 6$ and show that the integer in between the two limit points is either odd (not viable for this set of permutation) or when even, results in a non-integer value of b (Table 1). If it can be proven that $Epsilon(n) < Alpha_{min}(n)$ for all values of n than we have proven that cycling cannot occur with this particular set of iterations.

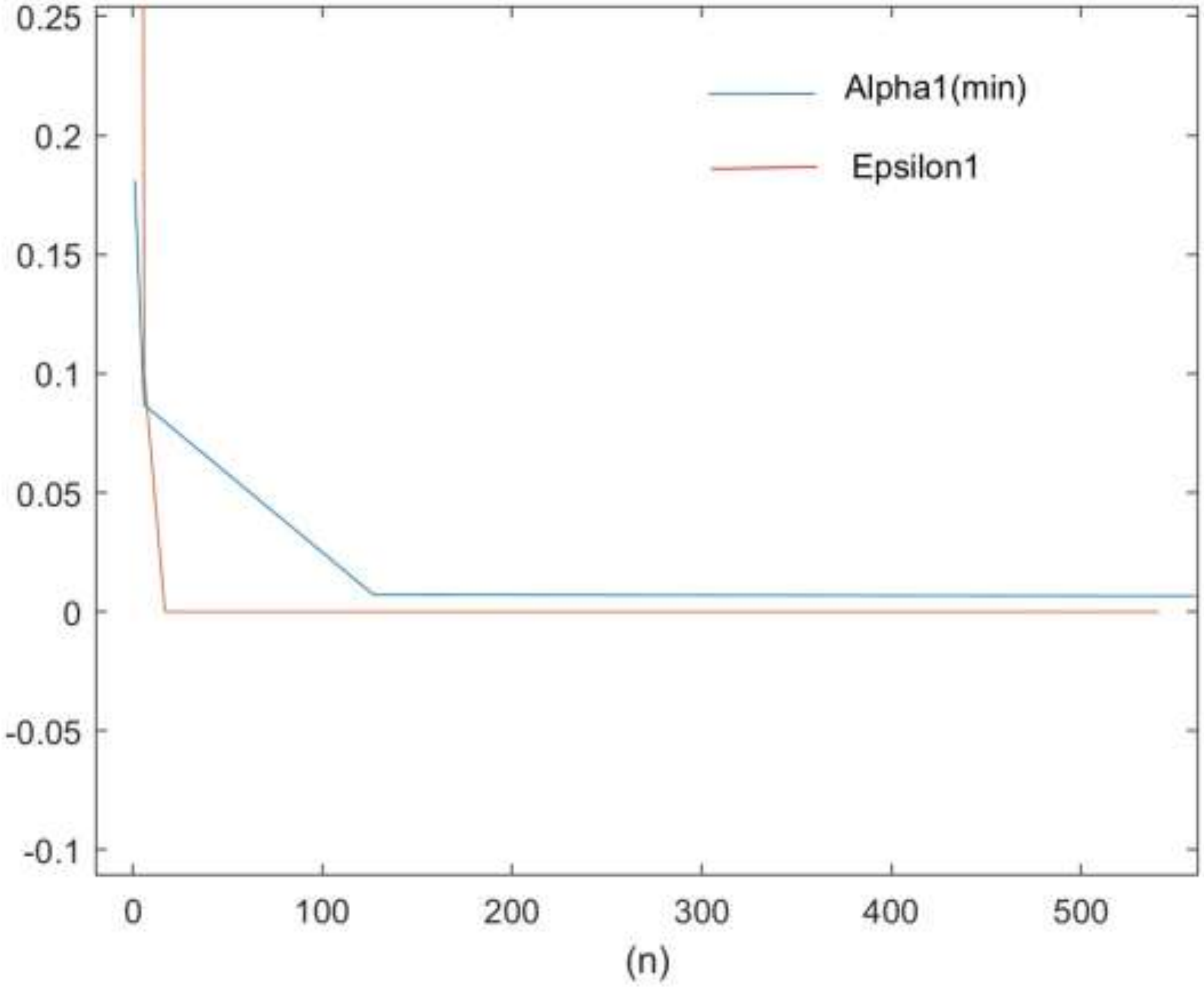

**Figure 2.** Alpha1(min) and Epsilon1 as a function of n. Epsilon1 is < Alpha1(min) indicates that there are no integers between the limit points of viable solutions for the cycling specified above. In this case it's true for values of n>6.

| n | $\frac{log_2 \frac{3^n}{2^n}}{(1-\frac{\varphi}{2})}$ | $\frac{log_2\left(\frac{3^n-1}{2^n-1}\right)}{(1-\frac{\varphi}{2})}$ | Difference | In Between Integer | $b = \frac{2^x - 3^{x/2}}{2^{x+n} - 3^{\frac{x}{2}+n}}$ |
|---|---|---|---|---|---|
| 1 | 2.82… | 4.82… | 2 | 4 | 1.4 |
| 2 | 5.64… | 6.83… | 1.19 | 6 | 2.84 |
| 3 | 8.46… | 9.12… | 0.67 | 9 | Integer not even |
| 6 | 16.91… | 17.01… | 0.09995 | 17 | Integer not even |

**Table 1.** Values of integers in between the two limit points when Epsilon1>Alpha1(min). When values of integers are even the resulting value of b is a non-integer. Otherwise, the values are odd which contradicts our assumption of this type of cycling.

**Second case**:

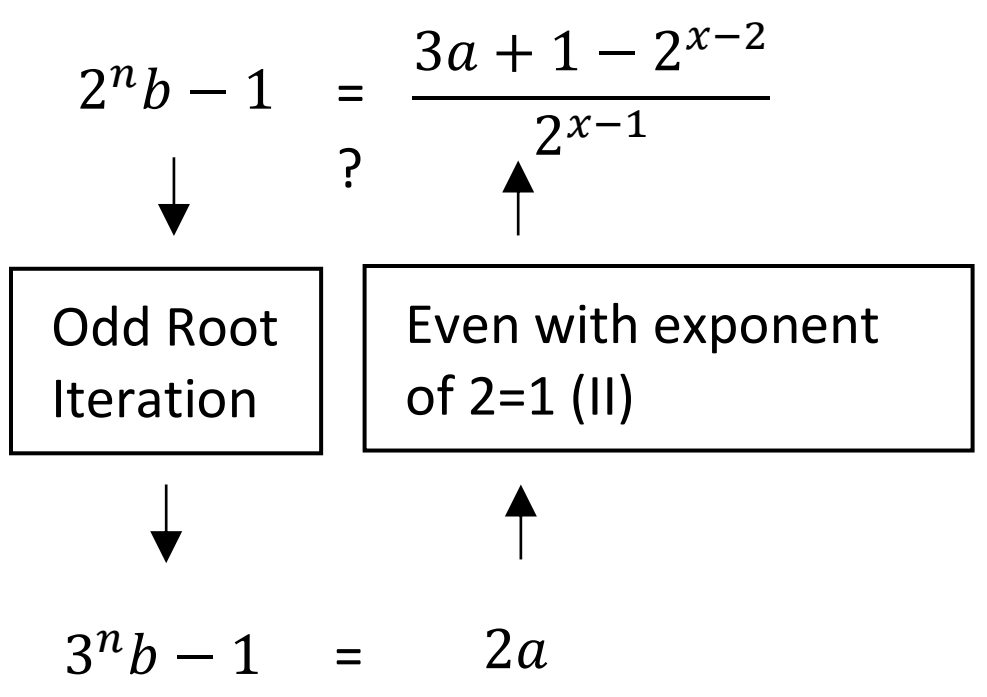

Again, I solve for b knowing that:

$$a = \frac{3^n b - 1}{2}$$

Therefore, if a cycle were to occur:

$$2^n b - 1 = \frac{3^{n+1} b - 3 + 2 - 2^{x-1}}{2^x}$$

Which simplifies to:

$$b = \frac{2^{x-1} - 1}{2^{n+x} - 3^{n+1}}$$

Next, I will Determine the values of x which will yield either b<0 or b<1 (i.e. non valid values of b)

b<0 when

$$2^{n+x} < 3^{n+1} \Longrightarrow 2^x < \frac{3^{n+1}}{2^n} \Rightarrow x < log_2\left(\frac{3^{n+1}}{2^n}\right).$$

b<1 when

$$2^{n+x} - 3^{n+1} > 2^{x-1} - 1 \Rightarrow 2^x > \frac{3^{n+1} - 1}{2^n - \frac{1}{2}}$$

Which simplifies to:

$$x > log_2\left(\frac{3^{n+1} - 1}{2^n - {}^1\!/_2}\right).$$

As previously, the two limit points will converge as n increase such that all values of x will lead to b<0 or b<1; leading us to the conclusion that this type of cycling cannot occur for large n. To illustrate that Epsilon is always less than Alpha (as defined previously) we design the following functions:

$$Alpha2_{min}(n) = \ min_1^n(1 - mod_1(log_2 \frac{3^{n+1}}{2^n}))$$

$$Epsilon2(n) = log_2\left(\frac{3^{n+1} - 1}{2^n - 1/2}\right) - log_2 \frac{3^{n+1}}{2^n}$$

In this case Epsilon2(n)<Alpha2(n) for all the calculated n's (Figure 3). It should be easier to prove in this case that Epsilon2(n)<Alpha2(n) for all n, since the equations are somewhat simpler.

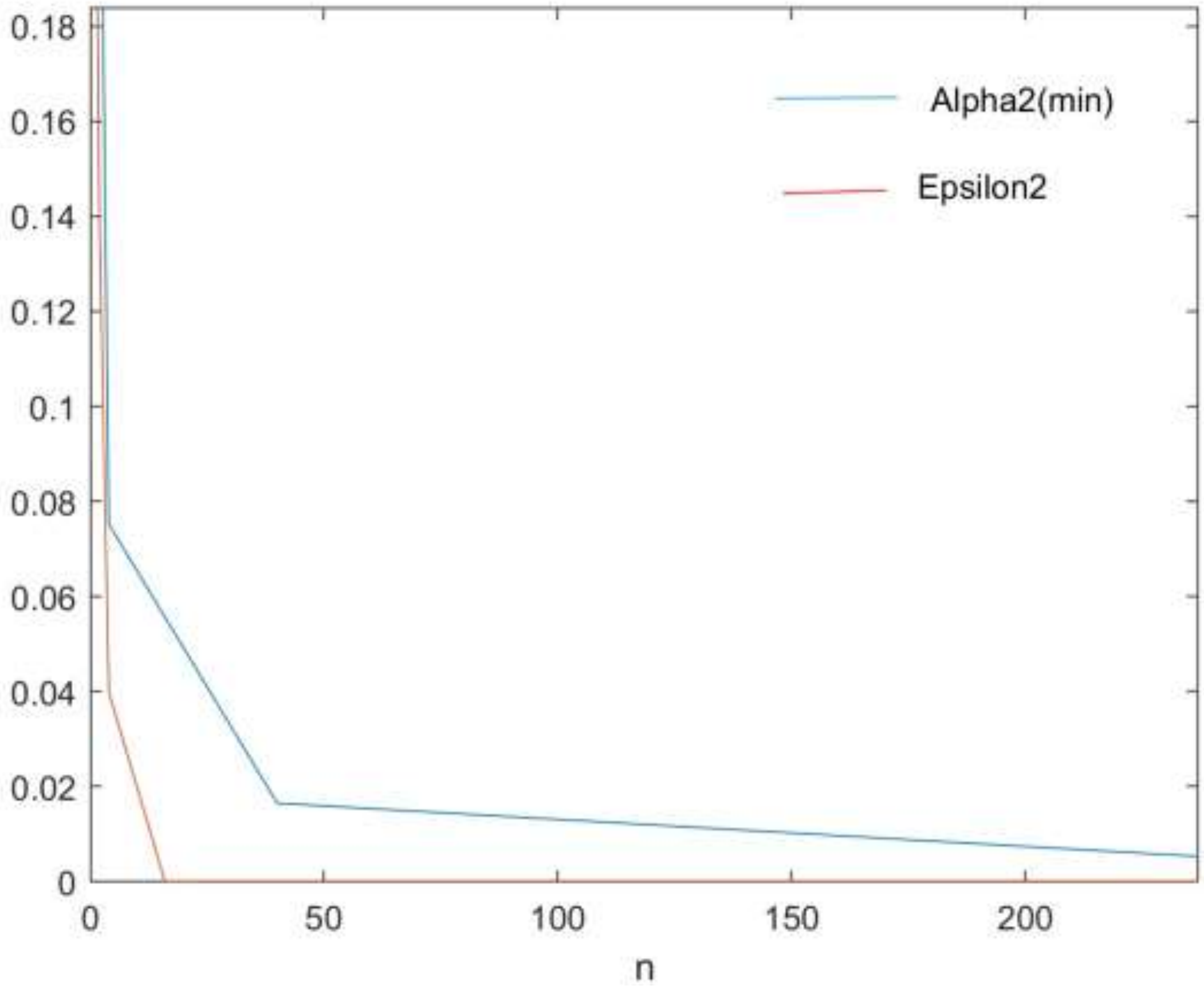

**Figure 3.** Alpha2(min) and Epsilon2 as a function of n. Epsilon2 is < Alpha3(min) indicates that there are no integers between the limit points of viable solutions for the cycling specified above. In this case it is true for all values of n.

**Third case:**

The last case we consider is when the sequence of iterations goes through the odd root iteration yielding an even root iteration of the form $2^x a$ but x is an odd exponent. The even root iteration will then yield the root iteration as diagramed below:

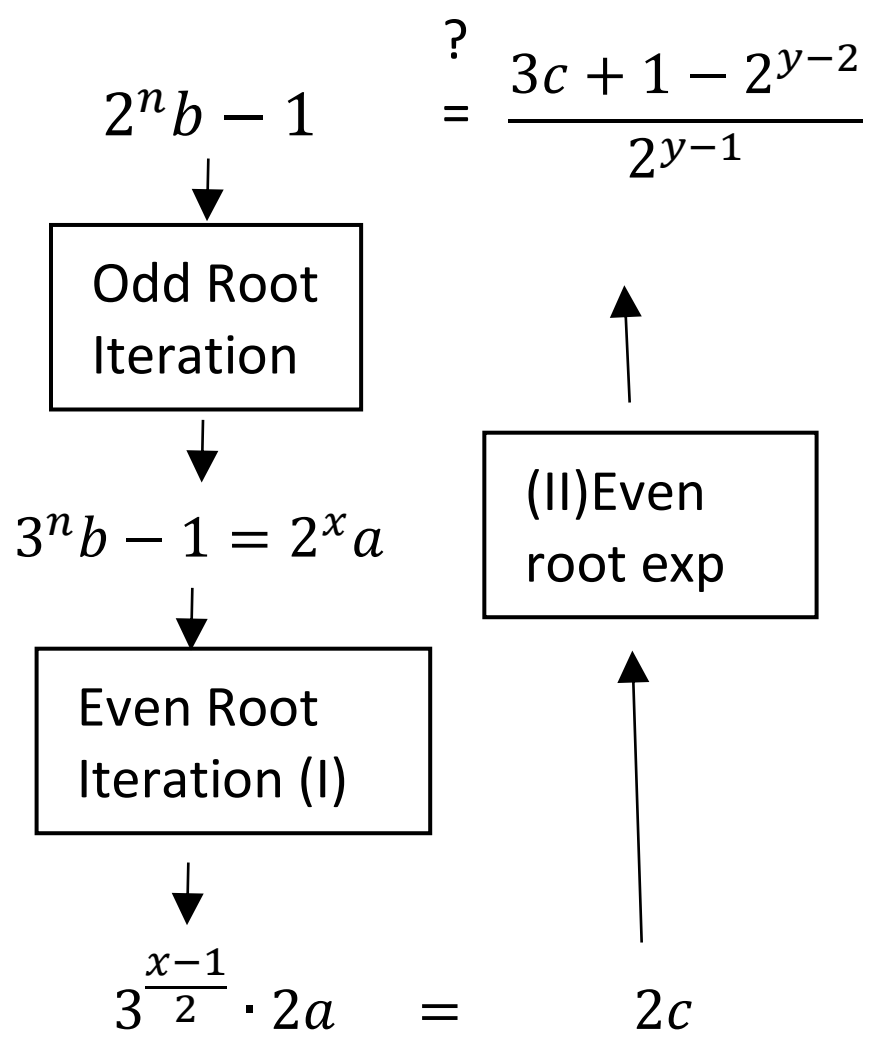

I will show, as in the other cases that for large n the limits converge leaving no plausible value of b; i.e. b<0 or <1 for all y's when n is large.

First, I describe c in terms of b.

$$c = 3^{\frac{x-1}{2}} \cdot a = 3^{\frac{x-1}{2}} \cdot \frac{3^n b - 1}{2^x} = \frac{3^{\frac{x-1}{2}+n} \cdot b - 3^{\frac{x-1}{2}}}{2^x}$$

Therefore:

$$\frac{3c+1-2^{y-2}}{2^{y-1}} = \frac{\frac{3^{\frac{x-1}{2}+n+1} \cdot b - 3^{\frac{x-1}{2}+1}}{2^x} + 1 - 2^{y-2}}{2^{y-1}}$$

Which reduces to:

$$= \frac{3^{\frac{x+1}{2}+n} \cdot b - 3^{\frac{x+1}{2}} \quad + 2^x - 2^{x+y-2}}{2^{x+y-1}}$$

If the above cycle exists

$$2^n b - 1 = \frac{3^{\frac{x+1}{2}+n} \cdot b - 3^{\frac{x+1}{2}} \quad + 2^x - 2^{x+y-2}}{2^{x+y-1}}$$

I will solve the above equation for b and examine what values of y yielding contradictory values of b (i.e. b<0 or b<1).

The above equation simplifies to:

$$b = \frac{2^x(1+2^{y-2}) - 3^{\frac{x+1}{2}}}{2^{x+y-1+n} - 3^{\frac{x+1}{2}+n}}$$

The values of y such that b<0 can be solved by setting the denominator of the above equation to a negative value (the numerator is positive).

$$2^{x+y-1+n} < 3^{\frac{x+1}{2}+n} \Rightarrow 2^y < \frac{3^{\frac{x+1}{2}+n}}{2^{x-1+n}}$$

If

$$y < log_2\left(\frac{3^{\frac{x+1}{2}+n}}{2^{x-1+n}}\right),$$

Then b<0.

Next, I ask the question: what are the values of y, which will give values of b less than 1? If b<1 then:

$$2^{x+y-1+n} - 3^{\frac{x+1}{2}+n} > 2^x(1+2^{y-2}) - 3^{\frac{x+1}{2}}$$

Solving the above inequality for y:

$$2^y > \frac{3^{\frac{x+1}{2}+n} - 3^{\frac{x+1}{2}}}{2^{x-1+n} - 2^{x-2}}$$

Which simplifies to

$$y > log_2\left(\frac{3^{\frac{x+1}{2}}(3^n-1)}{2^{x-2}(2^{n+1}-1)}\right) \ which\ converges\ to\ log_2\left(\frac{3^{\frac{x+1}{2}+n}}{2^{x-1+n}}\right)\ as\ n\ increases.$$

Hence, as in the other pathways, there are no plausible values of y as n increases. In this case it is more difficult to show the relationship between Alpha and Epsilon graphically as there are two independent variables in the equation: x and n.

I have shown that as the number of odd root iterations increases, the probability of the Collatz cycles shown here decreases. However, it must be proven that the Alpha and Epsilon functions behave as postulated here for every n. Interestingly, Eliahou (1993) has shown that any non-trivial cycle during the Collatz iteration must contain at least 17,087,915 elements. Such a large lower bound for possible cycles during the Collatz iteration and the findings presented here indicate that such large or larger cycles are non-existent. There are other possible cycles which must be considered. For example, there is a potential cycle where one has an odd-even-odd-even….. set of iterations, which must be shown to be unfeasible.

**<u>The Case for more complicated set of iterations.</u>**

The above possible cycles are simple and more complicated cycles must be considered. As mentioned above, one could have odd-even or even-odd cycles occurring before the iteration arrives at the initial number. Further the even iterations could be any combination of type I and type II even iterations (Figure 1). To address the possibility of cycling through more complicated set of iterations, an even iteration generic to both type I and II is designed to allow for any combination of even iterations.

**A generic even root iteration.**

Consider a number with an even root $= 2^j a$ where a is an odd number and $j \geq 1$. The number $(2^{x'+1}a + 1)$ can be iterated by the reduced Collatz iteration:

$$\frac{3(2^{j+1}a+1)+1}{2^x} = \frac{3 \cdot 2^{j+1}a+4}{2^x} = \frac{3 \cdot 2^{j-1}a+1}{2^{x-2}} = 2b+1$$

Solving for b the new root after one iteration:

$$\frac{3 \cdot 2^{j-1}a + 1 - 2^{x-2}}{2^{x-1}} = b$$

Notice that if $j \geq 2$ then $x = 2$ which reduces to the iteration given on Theorem 5

$$3 \cdot 2^{j-2}a$$

Examining multiple iterations with this generic form shows a pattern.

One Iteration:

$$\frac{3 \cdot 2^{j-1}a + (1 - 2^{x-2})}{2^{x-1}}$$

Two Iterations:

$$\frac{3^2 \cdot 2^{j-1}a + 3(1 - 2^{x_1-2}) + 2^{x_1-1}(1 - 2^{x_2-2})}{2^{x_1+x_2-2}}$$

Three Iterations:

$$\frac{3^3 \cdot 2^{j-1}a + 3^2(1 - 2^{x_1-2}) + 3 \cdot 2^{x_1-1}(1 - 2^{x_2-2}) + 2^{x_1+x_2-2}(1 - 2^{x_3-2})}{2^{x_1+x_2+x_3-3}}$$

Four Iterations:

$$\frac{3^4 \cdot 2^{j-1}a + 3^3(1 - 2^{x_1-2}) + 3^2 \cdot 2^{x_1-1}(1 - 2^{x_2-2}) + 3 \cdot 2^{x_1+x_2-2}(1 - 2^{x_3-2}) + 2^{x_1+x_2+x_3-3}(1 - 2^{x_4-2})}{2^{x_1+x_2+x_3+x_4-4}}$$

k Iterations:

$$\frac{3^k \cdot 2^{j-1}a + 3^{k-1}(1-2^{x_1-2}) + \cdots \ldots \ldots + 3 \cdot 2^{x_1+\cdots+x_{k-2}-(k-2)}(1-2^{x_{k-1}-2}) + 2^{x_1+x_2 \ldots + x_{k-1}-(k-1)}(1-2^{x_k-2})}{2^{x_1+x_2+x_3+\ldots+x_k-k}}$$

Or:

$$\frac{3^k \cdot 2^{j-1}a + \sum_{i=1}^{k} 3^{k-i} \cdot 2^{\sum_{l=0}^{i-1} x_l-(i-1)}(1-2^{x_i-2})}{2^{\sum_{i=1}^{k} x_i-k}}, where\ x_0 = 0$$

The numerator of the above equation can be thought of as having two parts: The Permanent part, which only changes by increasing the exponent of 3 after each iteration ($3^k \cdot 2^{j-1}a$) and the Elimination part, consisting of factors which are eliminated if $x_i = 2$ , which I make into a function to simplify the notation:

$$F(k,x) = F(k, x_1 \ldots\ldots x_k) = \sum_{i=1}^{k} 3^{k-i} \cdot 2^{\sum_{l=0}^{i-1} x_l-(i-1)}(1-2^{x_{i-2}})$$

Each of the above factor can be eliminated should $x_i = 2$.

Now I will examine the possibility of one cycle of odd root iterations followed by even root iterations using the generic iteration developed above.

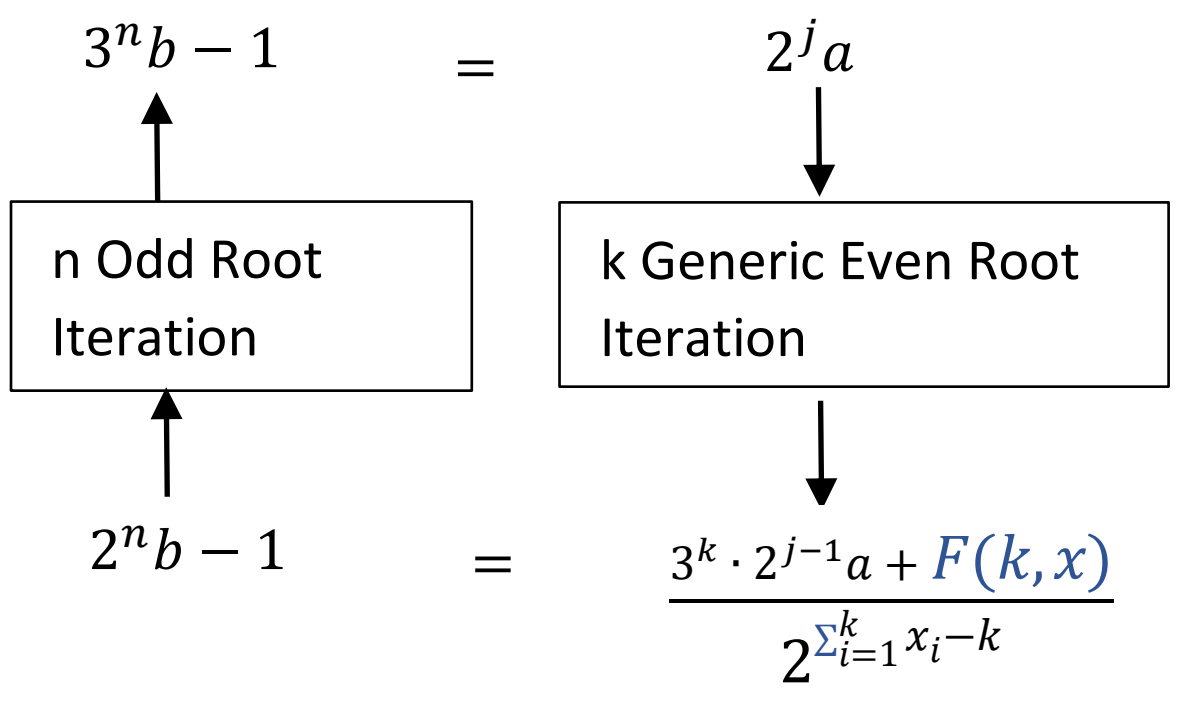

The critical equality is: $a = \frac{3^j b-1}{2^j}$

Substituting it into the following equation:

$$2^n b - 1 = \frac{3^k \cdot 2^{j-1}\frac{3^n b-1}{2^j} + F(k,x)}{2^{\sum_{i=1}^{k} x_i-k}}$$

$$= \frac{3^{k+n} \cdot 2^{j-1}b - 3^k \cdot 2^{j-1} + 2^j \cdot F(k,x)}{2^{\sum_{i=1}^{k} x_i-k+j}}$$

Therefore:

$$(2^{\sum_{i=1}^{k} x_i - k + j + n} - 3^{k+n} \cdot 2^{j-1})b = 2^{\sum_{i=1}^{k} x_i - k + j} - 3^k \cdot 2^{j-1} + 2^j \cdot F(k, x)$$

I examine the case where b<0.

$$2^{\sum_{i=1}^{k} x_i - k + j + n} < 3^{k+n} \cdot 2^{j-1}$$

$$\mathbf{2^{\sum_{i=1}^{k} x_i - k} < \frac{3^{k+n}}{2^{n+1}}}$$

For the case where b<1

$$2^{\sum_{i=1}^{k} x_i - k + j + n} - 3^{k+n} \cdot 2^{j-1} > 2^{\sum_{i=1}^{k} x_i - k + j} - 3^k \cdot 2^{j-1} + 2^j \cdot F(k, x)$$

Therefore

$$2^{\sum_{i=1}^{k} x_i - k + j + n} - 2^{\sum_{i=1}^{k} x_i - k + j} > 3^{k+n} \cdot 2^{j-1} - 3^k \cdot 2^{j-1} + 2^j \cdot F(k, x)$$

Factoring out $2^{\sum_{i=1}^{k} x_i - k + j}$

$$2^{\sum_{i=1}^{k} x_i - k + j}(2^n - 1) > 3^{k+n} \cdot 2^{j-1} - 3^k \cdot 2^{j-1} + 2^j \cdot F(k, x)$$

And

$$\mathbf{2^{\sum_{i=1}^{k} x_i - k} > \frac{3^k \cdot (3^n - 1)}{2(2^n - 1)} - \frac{F(k, x)}{2^n - 1}}$$

The values of the two inequalities for b<1 and b<0 will converge as n increases; showing that for large values of n there can be no viable cycling with odd and even root iterations.

Consider a more complex hypothetical cycle, by which two odd-even iteration cycles are needed to arrive at the same number. I abbreviate $n_i$ for the exponent, $a_i$ for the factor off the odd root in the ith cycle; $j_i$ for the exponent and $b_i$ the factor of the even root on the ith cycle; $x_i's$ are the exponents of the even iteration on the ith cycle; $k_i's$ are the number of even iterations in the ith cycle respectively.

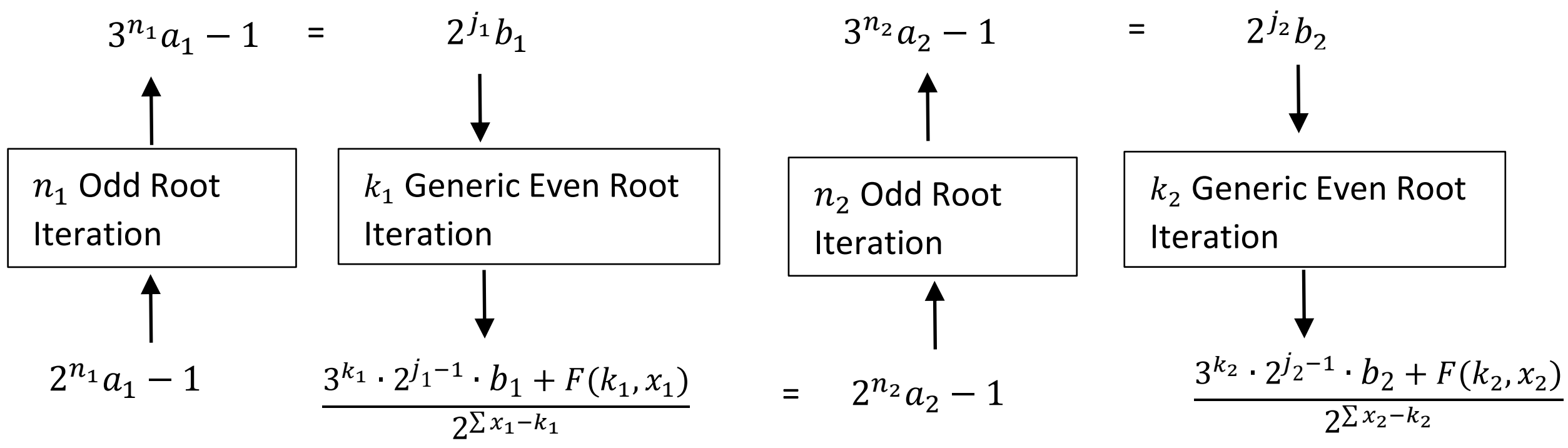

Critical equalities are:

$$b_1 = \frac{3^{n_1} a_1 - 1}{2^{j_1}}, \qquad a_2 = \frac{3^{k_1} \cdot 2^{j_1 - 1} \cdot b_1 + F(k_1, x_1) + 2^{\sum x_1 - k_1}}{2^{\sum x_1 - k_1 + n_2}}, \qquad b(2) = \frac{3^{n_2} a_2 - 1}{2^{j_2}}$$

After these two sets of iteration the sequence cycles back to $2^{n_1} a_1 - 1$, therefore:

$$2^{n_1} a_1 - 1 = \frac{3^{k_2} \cdot 2^{j_2 - 1} \cdot b_2 + F(k_2, x_2)}{2^{\sum x_2 - k_2}}$$

Substituting the value of $b_2$.

$$2^{n_1} a_1 - 1 = \frac{3^{k_2} \cdot 2^{j_2 - 1} \cdot \frac{3^{n_2} a_2 - 1}{2^{j_2}} + F(k_2, x_2)}{2^{\sum x_2 - k_2}} = \frac{3^{k_2 + n_2} \cdot 2^{j_2 - 1} a_2 - 3^{k_2} \cdot 2^{j_2 - 1} + 2^{j_2} F(k_2, x_2)}{2^{\sum x_2 - k_2 + j_2}}$$

Substituting the value of $a_2$.

$$= \frac{3^{k_2 + n_2} \cdot 2^{j_2 - 1} \frac{3^{k_1} \cdot 2^{j_1 - 1} \cdot b_1 + F(k_1, x_1) + 2^{\sum x_1 - k_1}}{2^{\sum x_1 - k_1 + n_2}} - 3^{k_2} \cdot 2^{j_2 - 1} + 2^{j_2} F(k_2, x_2)}{2^{\sum x_2 - k_2 + j_2}}$$

$$= \frac{3^{k_2 + n_2 + k_1} \cdot 2^{j_2 + j_1 - 2} b_1 + 3^{k_2 + n_2} \cdot 2^{j_2 - 1} F(k_1, x_1) + 3^{k_2 + n_2} \cdot 2^{j_2 - 1 + \sum x_1 - k_1} - 3^{k_2} \cdot 2^{j_2 - 1 + \sum x_1 - k_1 + n_2} + 2^{j_2 + \sum x_1 - k_1 + n_2} F(k_2, x_2)}{2^{\sum x_2 - k_2 + j_2 + \sum x_1 - k_1 + n_2}}$$

Substituting the value $b_1$ of in terms of $a_1$.

$$= \frac{3^{k_2 + n_2 + k_1} \cdot 2^{j_2 + j_1 - 2} \frac{3^{n_1} a_1 - 1}{2^{j_1}} + 3^{k_2 + n_2} \cdot 2^{j_2 - 1} F(k_1, x_1) + 3^{k_2 + n_2} \cdot 2^{j_2 - 1 + \sum x_1 - k_1} - 3^{k_2} \cdot 2^{j_2 - 1 + \sum x_1 - k_1 + n_2} + 2^{j_2 + \sum x_1 - k_1 + n_2} F(k_2, x_2)}{2^{\sum x_2 - k_2 + j_2 + \sum x_1 - k_1 + n_2}}$$

$$= \frac{3^{k_2+n_2+k_1+n_1} \cdot 2^{j_2+j_1-2} a_1 - 3^{k_2+n_2+k_1} 2^{j_2+j_1-2} + 3^{k_2+n_2} \cdot 2^{j_2+j_1-1} F(k_1, x_1) + 3^{k_2+n_2} \cdot 2^{j_2-1+\sum x_1-k_1+j_1} - 3^{k_2} \cdot 2^{j_2-1+\sum x_1-k_1+n_2+j_1} + 2^{j_2+\sum x_1-k_1+n_2+j_1} F(k_2, x_2)}{2^{\sum x_2-k_2+j_2+\sum x_1-k_1+n_2+j_1}}$$

Therefore:

$$\left(2^{\sum x_2-k_2+j_2+\sum x_1-k_1+n_2+j_1+n_1} - 3^{k_2+n_2+k_1+n_1} \cdot 2^{j_2+j_1-2}\right) a_1 =$$

$$2^{\sum x_2-k_2+j_2+\sum x_1-k_1+n_2+j_1} - 3^{k_2+n_2+k_1} 2^{j_2+j_1-2} + 3^{k_2+n_2} \cdot 2^{j_2+j_1-1} F(k_1, x_1) + 3^{k_2+n_2} \cdot 2^{j_2-1+\sum x_1-k_1+j_1} - 3^{k_2} \cdot 2^{j_2-1+\sum x_1-k_1+n_2+j_1} + 2^{j_2+\sum x_1-k_1+n_2+j_1} F(k_2, x_2)$$

Now examine the case where *a<0*

$$2^{\sum x_2-k_2+j_2+\sum x_1-k_1+n_2+j_1+n_1} < 3^{k_2+n_2+k_1+n_1} \cdot 2^{j_2+j_1-2}$$

$$\mathbf{2^{\sum x_2-k_2+\sum x_1-k_1+n_2} < \frac{3^{k_2+n_2+k_1+n_1}}{2^{n_1+2}}}$$

Now examine the case where *a<1*

$$2^{\sum x_2-k_2+j_2+\sum x_1-k_1+n_2+j_1+n_1} - 3^{k_2+n_2+k_1+n_1} \cdot 2^{j_2+j_1-2} >$$

$$2^{\sum x_2-k_2+j_2+\sum x_1-k_1+n_2+j_1} - 3^{k_2+n_2+k_1} 2^{j_2+j_1-2} + 3^{k_2+n_2} \cdot 2^{j_2+j_1-1} F(k_1, x_1) + 3^{k_2+n_2} \cdot 2^{j_2-1+\sum x_1-k_1+j_1} - 3^{k_2} \cdot 2^{j_2-1+\sum x_1-k_1+n_2+j_1} + 2^{j_2+\sum x_1-k_1+n_2+j_1} F(k_2, x_2)$$

$$2^{\sum x_2-k_2+j_2+\sum x_1-k_1+n_2+j_1+n_1} - 2^{\sum x_2-k_2+j_2+\sum x_1-k_1+n_2+j_1} >$$

$$3^{k_2+n_2+k_1+n_1} \cdot 2^{j_2+j_1-2} - 3^{k_2+n_2+k_1} 2^{j_2+j_1-2} + 3^{k_2+n_2} \cdot 2^{j_2+j_1-1} F(k_1, x_1) + 3^{k_2+n_2} \cdot 2^{j_2-1+\sum x_1-k_1+j_1} - 3^{k_2} \cdot 2^{j_2-1+\sum x_1-k_1+n_2+j_1} + 2^{j_2+\sum x_1-k_1+n_2+j_1} F(k_2, x_2)$$

$$2^{\sum x_2 - k_2 + j_2 + \sum x_1 - k_1 + n_2 + j_1}(2^{n_1} - 1) >$$

$$3^{k_2+n_2+k_1} 2^{j_2+j_1-2}(3^{n_1} - 1) + 3^{k_2+n_2} \cdot 2^{j_2+j_1-1} F(k_1, x_1) + 3^{k_2+n_2} \cdot 2^{j_2-1+\sum x_1 - k_1 + j_1} - 3^{k_2} \cdot 2^{j_2-1+\sum x_1 - k_1 + n_2 + j_1} + 2^{j_2+\sum x_1 - k_1 + n_2 + j_1} F(k_2, x_2)$$

$$2^{\sum x_2 - k_2 + \sum x_1 - k_1 + n_2} > \frac{3^{k_2+n_2+k_1}(3^{n_1} - 1)}{2^2(2^{n_1} - 1)} + \frac{3^{k_2+n_2} \cdot 2^{-1} F(k_1, x_1) + 3^{k_2+n_2} \cdot 2^{\sum x_1 - k_1 - 1} - 3^{k_2} \cdot 2^{\sum x_1 - k_1 + n_2 - 1} + 2^{\sum x_1 - k_1 + n_2} F(k_2, x_2)}{(2^{n_1} - 1)}$$

$$\mathbf{2^{\sum x_2 - k_2 + \sum x_1 - k_1 + n_2} > \frac{3^{k_2+n_2+k_1}(3^{n_1} - 1)}{2^2(2^{n_1} - 1)} + \frac{N_4}{2^{n_1} - 1}}$$

It's easy to see that the $2^{\sum x_2 - k_2 + \sum x_1 - k_1 + n_2}$ converge to the same value for $a_1 < 0$ and $a_1 < 1$ as $n_1$ increases. Therefore, all values of $2^{\sum x_2 - k_2 + \sum x_1 - k_1 + n_2}$ will either yield $a_1 < 0$ or $a_1 < 1$ , when $n_1$ is large, contradicting the assumption that $a_1$ is an integer.

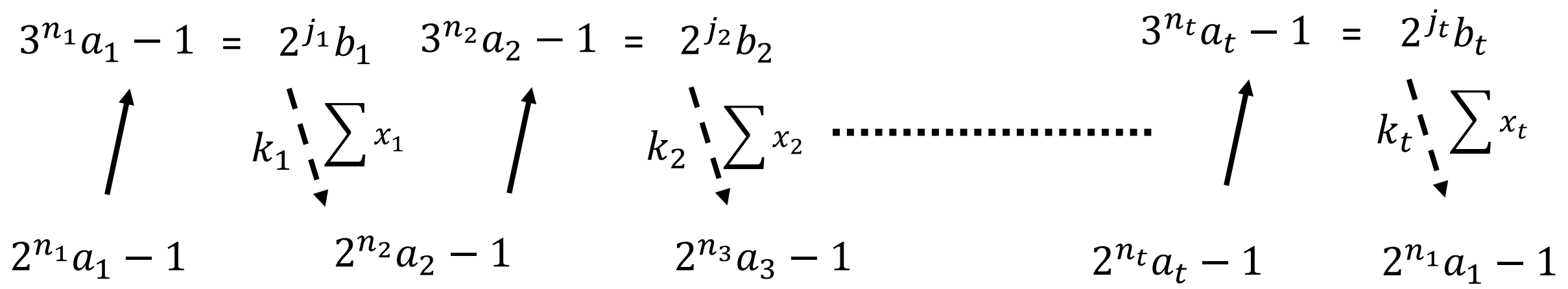

**Figure 5.** t cycles of odd – even iterations. $n_i$ represents the exponent of the multiplicative factor of odd roots, $a_i$ represents the other term in the multiplicative factor of odd roots, $j_1$ represents the exponent of the even root before interpolations with the odd term $b_i$, $k_i$ represents the number of even root iterations on the ith cycle and $\sum x_i$ represents the sum of all the exponents in the even iteration on the ith cycle.

It's easy to see the extrapolation of the limiting inequalities for t cycles when $a_i < 0$

$$\mathbf{2}^{\sum x_t - k_t + \sum x_{t-1} - k_{t-1} + n_t \dots\dots\dots\dots \sum x_1 - k_1 + n_2} < \frac{\mathbf{3}^{k_t + n_t \dots\dots\dots\dots + k_1 + n_1}}{\mathbf{2}^{n_1 + t}}$$

And for the case when $a_i < 1$

$$\mathbf{2}^{\sum x_t - k_t + \sum x_{t-1} - k_{t-1} + n_t \dots\dots\dots\dots \sum x_1 - k_1 + n_2} > \frac{\mathbf{3}^{k_t + n_t \dots\dots\dots + k_1}(\mathbf{3}^{n_1} - \mathbf{1})}{\mathbf{2}^t(\mathbf{2}^{n_1} - \mathbf{1})} + \frac{N_t}{\mathbf{2}^{n_1} - \mathbf{1}}$$

These two inequalities will also converge to the same value as $n_1$increases leading to the conclusion that regardless of the number of cycles and types of iterations there is no cycling as long as $n_1$ is sufficiently large to render Epsilon(t) less than Alpha(t) as defined previously. It's assumed above that t cycles of odd-even iterations will lead to the initial value of the odd root. Therefore, the prerequisite that $n_1$is sufficiently large to render it impossible for the odd even iterations to cycle can be broadened to say that if any $n_i$ is sufficiently large relative to the other exponents, than there is no cycling possible. In other words, if the multiplicative factor of any of the odd roots in the hypothetical cycle above is of a sufficiently high power of two, then no cycling can occur.

**Starting with an even root**.

Previous only cycles that initiated with a odd rooted number were considered.  Here I consider a cycle that starts with an even rooted number goes through a set of even and odd interpolations to arrive at the same even rooted number:

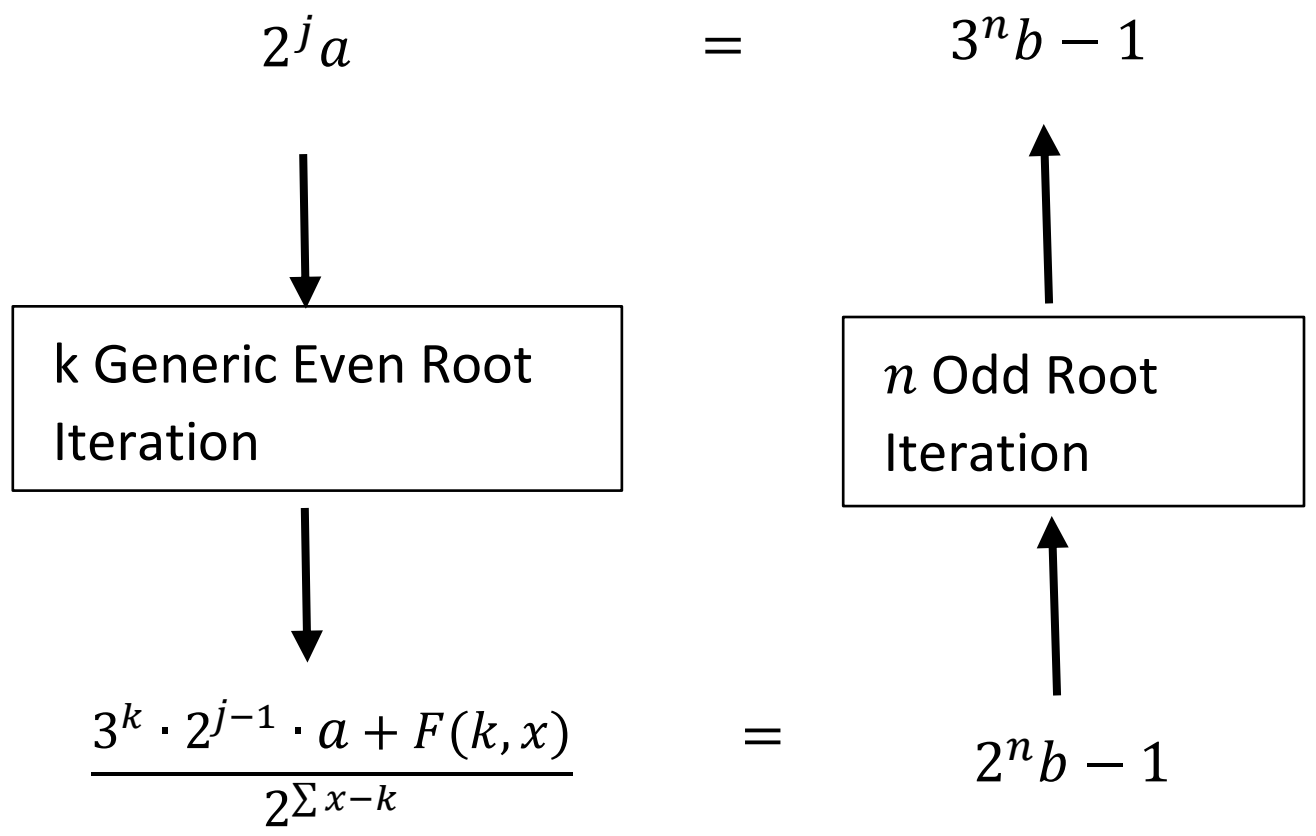

Critical equality is:

$$b = \frac{3^k \cdot 2^{j-1} \cdot a + F(k,x) + 2^{\sum x-k}}{2^{\sum x-k+n}}$$

Therefore:

$$2^j a = 3^n b - 1 = 3^n \frac{3^k \cdot 2^{j-1} \cdot a + F(k,x) + 2^{\sum x-k}}{2^{\sum x-k+n}} - 1$$

$$= \frac{3^{k+n} \cdot 2^{j-1} \cdot a + 3^n F(k,x) + 3^n 2^{\sum x-k} - 2^{\sum x-k+n}}{2^{\sum x-k+n}}$$

$$2^{j+\sum x-k+n} a = 3^{k+n} \cdot 2^{j-1} \cdot a + 3^n F(k,x) + 3^n 2^{\sum x-k} - 2^{\sum x-k+n}$$

$$(2^{j+\sum x-k+n} - 3^{k+n} \cdot 2^{j-1}) a = 3^n F(k,x) + 3^n 2^{\sum x-k} - 2^{\sum x-k+n}$$

For the case of $a < 0$

$$2^{j+\sum x-k+n} < 3^{k+n} \cdot 2^{j-1}$$

And

$$\mathbf{2^{\sum x-k+n} < \frac{3^{k+n}}{2}}$$

For the case where  $a < 1$

$$2^{j+\Sigma x-k+n} - 3^{k+n} \cdot 2^{j-1} > 3^n F(k,x) + 3^n 2^{\Sigma x-k} - 2^{\Sigma x-k+n}$$

$$2^{j+\Sigma x-k+n} + 2^{\Sigma x-k+n} > 3^{k+n} \cdot 2^{j-1} + 3^n F(k,x) + 3^n 2^{\Sigma x-k}$$

$$2^{\Sigma x-k+n}(2^j + 1) > 3^{k+n} \cdot 2^{j-1} + 3^n F(k,x) + 3^n 2^{\Sigma x-k}$$

$$\mathbf{2}^{\Sigma x-k+n} > \frac{\mathbf{3}^{k+n} \cdot \mathbf{2}^{j-1}}{(\mathbf{2}^j + \mathbf{1})} + \frac{\mathbf{3}^n F(k,x) + \mathbf{3}^n \mathbf{2}^{\Sigma x-k}}{(\mathbf{2}^j + \mathbf{1})}$$

Therefore as $j$ increases these two limit points will converge and no viable exponent of 2 during iterations will produce a positive integer $a$. The interpolations starting with an even rooted number can be expanded to t cycles as I did with odd rooted numbers. Given the observations for both even and odd rooted numbers above one can make the generalization that: given a set of cycles if one of the root multiplicative factor (be it odd or even) is a high power of 2 ($j$ $or$ $n$) no cycling is possible.